\documentclass[12pt,oneside,reqno]{amsart}

\usepackage{amsmath}
\usepackage{amssymb}
\usepackage{amsthm}

\newtheorem{theorem}{Theorem}
\newtheorem{corollary}{Corollary}
\newtheorem{lemma}{Lemma}

\newtheorem*{conjecture}{Conjecture}

\theoremstyle{definition}
\newtheorem{definition}{Definition}
\newtheorem{remark}{Remark}
\newtheorem{example}{Example}

\newcommand{\N}{\mathbb N}
\newcommand{\R}{\mathbb R}

\newcommand{\Sim}{\operatorname{Sim}}
\newcommand{\Isom}{\operatorname{Isom}}
\newcommand{\Hd}{\operatorname{Hd}}
\newcommand{\Id}{\operatorname{Id}}
\newcommand{\const}{\operatorname{const}}
\newcommand{\hb}{\mathcal{HB}}

\begin{document}

\large

\title{Semilinear metric semilattices on $\R$-trees}
\author{P.D. Andreev}
\address{Lomonosov Pomor State University, Arkhangelsk, Russia}
\thanks{Supported by RFBR, grant 04-01-00315-a and program of Department of education and science of Russia "Development of the scientific potential of high school", code of project 335.}
\email{pdandreev@mail.ru}
\maketitle
\begin{abstract}
We introduce the notion of metric semilattice on the metric space and prove the criterion of $\R$-tree as connected geodesic metric space $X$ admitting the partial order, such that $X$ is semilinear metric semilattice. Also we state the homeomorphism between topological space of orders defining upper semilinear metric $\vee$-semilattices on locally compact complete $\R$-tree $X$ and its metric compactification $\overline{X}_m$. As an application we construct the example of locally complete non-homogeneous similarity-homogeneous space showing essentiality of the condition of locally compactness in V.N. Berestovski\v\i's conjecture on the structure of such spaces. Constructed metric space is $\R$-tree, where every point is a branching point. It is the metric fibration but is not topological product with factor $\R$ and does not satisfy the Berestovski\v\i's conjecture.
\end{abstract}

\section{Introduction}

The notion of $\R$-tree is the generalization of symplicial tree and it is included into more general family of so called $\Lambda$-trees. The geodesic metric space $X$ is called $\R$-tree, if for every triangle any its side is contained in the union of another two. Here we introduce the notion of metric semilattice and prove following criterion:

\begin{theorem}\label{crit}
Let $X$ be a geodesic metric space. If $X$ is an $\R$-tree, then for every point $o \in X$ there exists unique partial order $\preceq$ on $X$ such that the pair $(X,\preceq)$ is upper-semilinear $\vee$-semilattice with root $o$. Every nonempty subset $A \subset X$ has its supremum in this order. Conversely, if $X$ admits partial order, such that $X$ becomes upper-semilinear metric semilattice, where directions of semilattice and semilinearity coincide, then $X$ is an $\R$-tree.
\end{theorem}

In the paragraph \ref{spaceoforders} we study the set $\mathcal O_+(X)$ of partial orders on the complete locally compact $\R$-tree $X$, which define on $X$ upper-semilinear metric  $\vee$-semilattices. We introduce the topology in $\mathcal O_+(X)$. In subspace $\mathcal O_+^r(X) \subset \mathcal O_+(X)$ consisting of rooted orders this topology is generated by Hausdorff metric on the family $\mathcal C(X\times X)$ of closed subsets of metric square $X \times X$. The continuation of the topology to the entire $\mathcal O_+(X)$ is constructed within the base of neighbourhoods of non-rooted orders. The theorem is proved.

\begin{theorem}
The metric space $\mathcal O_+^r(X)$ is isometric to $X$, and the topological space $\mathcal O_+(X)$ is homeomorphic to the metric compactification $\overline{X}_m$ of $X$.
\end{theorem}

As an application of the theorem \ref{crit} the example, showing essentiality of the condition of locally compactness in following conjecture formulated by V.N. Berestovski\v\i\  in \cite{Be} is built in paragraph \ref{example}.

\begin{conjecture}
Every locally compact similarity-homogeneous non-homogeneous metric space with inner metric  $(X, \rho)$ is homeomorphic to topological product $F \times \R_+$, where $F$ is arbitrary level set of the function of radius of completeness on $X$. Topological group $\Sim(X)$ of similarities of $X$ is homeomorphic to direct topological product $\Isom(X) \times \R_+$, where $\Isom(X) \subset \Sim(X)$ is the subgroup of isometries of $X$.
\end{conjecture}

Here $\R_+$ denotes the set of positive real numbers. The metric space $X$ constructed here is $\R$-tree and satisfies all conditions of the conjecture except locally compactness. It is not homeomorphic to topological product $F \times \R_+$ but is metric fibration in the sense of definition of the paper \cite{VGE}.  Every point of the space is branching point. The group of similarities $\Sim(X)$ splits by the exact sequence
\begin{equation}\label{split}
0 \to \Isom(X) \to \Sim(X) \to \R_+ \to 0,
\end{equation}
but is not homeomorphic to topological group $\Isom(X) \times \R_+$.

\section{Preliminaries}

Let $(X, \rho)$ be a metric space. \textit{The minimizer} in $X$ is by definition an image in the map $\gamma : (\alpha, \beta) \to X$, (\textit{natural parameterization}  of the minimizer) of real interval $(\alpha, \beta)$, such that $\rho(\gamma(s)\gamma(t)) = |s-t|$ for all $s,t \in (\alpha, \beta)$. If $(\alpha, \beta)$ in the definition above is the real segment, we say that $\gamma$ is \textit{segment} in $X$ connecting $x = \gamma(\alpha)$ with $y = \gamma(\beta)$, if $\alpha \in \R$ and $\beta = +\infty$, then $\gamma$ is \textit{ray}, and if $(\alpha, \beta) = \R$, then $\gamma$ is a \textit{straight line}. The natural parameterization of the minimizer is defined up to addition of constant.

The metric space $X$ is called \textit{geodesic} if every two its points $x, y \in X$ can be connected by segment. In general such segment is not unique. The space $X$ is called \textit{locally complete}, if for every point $x \in X$ there exists such a number $c(X) > 0$ that all balls centered in $x$ with radii less than $c(x)$ are complete. The maximal number $c(x)$ with this property is called \textit{the completeness radius} in $x$. The completeness radius is continuous function of point $x \in X$.

The map $f : X \to Y$ of metric space $(X,\rho_X)$ to metric space $(Y, \rho_Y)$ is called \textit{similarity} with coefficient $k > 0$, if for all $x, y \in X$ the equality $\rho_Y(f(x), f(y)) = k \cdot \rho_X(x,y)$ holds. When $k=1$ the similarity is \textit{isometric map}. The similarity map of $X$ onto itself is the similarity of $X$ and isometric map of $X$ onto itself is \textit{isometry} of $X$.

The geodesic metric space $X$ is called \textit{$\R$-tree} if every two its points are connected by unique segment and for every triple of points $x, y, z \in X$ the segment $[xy]$ is contained in the union of segments $[xz] \cup [yz]$. Analogous inclusion is true for segments $[xz]$ and $[yz]$ as well. The survey of the theory of $\R$-trees can be found in the article \cite{Bes}.

For complete locally compact metric space $(X,\rho)$ its metric compactification $\overline{X}_m$  is defined in \cite{Gr}. It admits following description. Consider Kuratowskii's embedding of the space $X$ into the space $С(X, \R)$ of continuous functions equipped with compact-open topology. Every point $x \in X$ maps to its distance function $d_x$:
\[d_x(y) = \rho(x,y) - \rho(o,x),\]
where $o \in X$ is marked point. Changing of marked point $o$ leads to addition of constants to all distance functions, so one can continue the embedding to the embedding of $X$ into factorspace $C^\ast(X, \R) = C(X, \R)/ \{consts\}$ of vector space $C(X, \R)$ by its subspace of constants. $X$ is identified with its image in $C^\ast(X,\R)$. \textit{The metric compactification} of the space $X$ is its closure in $C^\ast(X, \R)$. The boundary $\partial_mX = \overline{X}_m \setminus X$ is called {\em metric boundary}. Limiting functions in $\partial_mX = \overline{X} \setminus X$ are \textit{horofunctions}.

The equivalent definition of metric compactification as state space of unital commutative $C^\ast$-algebra generated by constants, functions vanishing at infinity and differences of distance functions is given in \cite{WW}.

The space $X$ is called \textit{homogeneous} (correspondingly, \textit{similarity homogeneous}), if its group of isometries (correspondingly similarities) acts on $X$ transitively. In \cite{Be} the metric structure of locally complete similarity homogeneous non-homogeneous metric spaces with inner metric is studied.  It is shown that every such space is conformally equivalent to complete homogeneous space with inner metric. The function of completeness radius in this case is a submetry of the space $X$ onto $\R_+$, i.e. it maps arbitrary ball in $X$ onto the ball of the same radius in $\R_+$.

Two subsets $F_1, F_2$ of the metric space $X$ are called \textit{equidistant} if for any point $x_i \in F_i, \ i = \overline{1,2}$ there exists a point $x_j \in F_j, \ j \ne i$ for which the distance $\rho(x_i x_j)$ is equal to the distance between $F_1$ and $F_2$. \textit{The metric fibration} $\mathcal F$ of the space $X$ is its subdivision into the family of mutually isometric relatively metric induced by $\rho$ equidistant closed sets. The factor-set $M/\mathcal F$ inherits natural metric $\rho(F_1, F_2), F_1, F_2 \in \mathcal F$, for which the factorization map $p: M \to M/ \mathcal F$ is a submetry.

For arbitrary metric space $(X, \rho)$ there is defined so called \textit{Hausdorff metric}  $\Hd$ on the family $\mathcal C(X)$ of closed sets in $X$ which can be infinite. The Hausdorff distance between closed sets $V,W \subset X$ is the value
\[\Hd(V,W) = \inf\{\varepsilon > 0\ |\ V \subset \mathcal N_\varepsilon(W),\ W \subset \mathcal N_\varepsilon(V)\},\]
where $\mathcal N_\varepsilon(P)$ denotes $\varepsilon$-neighbourhood of the set $P\subset X$, i.e.
\[\mathcal N_\varepsilon(P) = \{y \in X\ |\ \exists x \in P; \ \rho(x,y) < \varepsilon\}.\]

\section{The $\R$-tree criterion}

In this paragraph we give the definition of metric semilattice and prove the theorem \ref{crit}. Main facts of the theory of posets and lattices can be found in \cite{Bi}.

\begin{definition}\label{dist}
Let the partial order $\preceq$ and corresponding strong order $\prec$ are given on the metric space $(X, \rho)$. Suppose the poset $(X, \preceq)$ is $\vee$-semilattice ($\wedge$-semilattice), i. e. for every two points $x,y \in X$ there exists their supremum (correspondingly infimum) $x \vee y$ (correspondingly $x\wedge y$). The triple $(X, \rho, \preceq)$ is called \textit{metric $\vee$-semilattice} (correspondingly, \textit{metric  $\wedge$-semilattice}), if 
\begin{enumerate}
\item for every $x, y, z \in X$ from $x \preceq z \preceq y$ it follows $\rho(x,z) + \rho(z,y) = \rho(x,y)$, and 
\item
for every two points $x, y \in X$ equality holds:
\[\rho(x,y) = \rho(x, x\vee y) + \rho(x \vee y, y)\]
(correspondingly,
\[\rho(x,y) = \rho(x, x\wedge y) + \rho(x \wedge y, y)).\]
\end{enumerate}
\end{definition}

Since, because of duality principle, every proposition about $\vee$-semilattices is automatically valid for $\wedge$-semilattices, we will deal only with $\vee$-semilattices.

Recall that the partial order $\preceq$ on the set $X$ is called \textit{upper (low) semilinear}, if for every point $x \in X$ its upper (correspondingly low) cone
\[U_x = \{y \in X\ | \ x \preceq y\}\]
(correspondingly,
\[L_x = \{z \in X\ |\ z \preceq x\})\]
is linearly ordered.

\textit{Proof of the theorem \ref{crit}. Existence of the order}. Let $X$ be $\R$-tree. Take a point $o \in X$ and define relation $\preceq_{(o)}$ by condition: for points $x, y \in X$ put $y \preceq_{(o)} x$ when and only when $x \in [oy]$. The relation $\preceq_{(o)}$ is transitive because if $x \in [oy]$ and $y \in [oz]$, then by definition of $\R$-tree $[oy] \subset [oz]$ and $x \in [oz]$. Reflexivity and antisymmetry of $\preceq_{(o)}$ are evident. Hence the relation $\preceq_{(o)}$ is partial order.

The order $\preceq_{(o)}$ is upper-semilinear: if $x \preceq_{(o)} y$ and $x \preceq_{(o)} z$, then $y, z \in [ox]$ and $y \in [oz]$ or $z \in [ox]$. Hence all cones $U_x$ are linearly ordered.

Consider arbitrary subset $B \subset X$ and set
\[Z(B) = \bigcap_{b \in B}[ob].\] 
The set $Z(B)$ is nonempty because of $o \in Z(B)$. Moreover, $Z(B)$ is one-point set or is homeomorphic to segment. Really:  $Z(B)$ is contained in every segment $[ob]$ and linearly connected, so consists of unique point $o$, or is homeomorphic to real interval. In the latter we consider arbitrary point $b \in B$ and the point $x$ of segment $[ob]$ which is the most distant from $o$ limit point of the interval $Z(B)$. If $b' \in B$ is the point, different from $b$, then $[ox] \subset [ob] \cap [ob']$. Hence $x \in Z(B)$ and $Z(B) = [ox]$. The supremum $\bigvee_{(o)}B$ of the set $B$ is the point $o$ when $Z(B) = \{o\}$ or $x$ when $Z(B) = [ox]$. Supremum of points $a,b \in X$ is denoted $a \vee_{(o)} b$.

The order $\preceq_{(o)}$ is metric. Really, the proposition (1) in the definition \ref{dist} is true automatically. Moreover, for $a, b \in X$
\[Z(\{a,b\})= \left\{
	\begin{array}{ll}
	~\{o\},& \mbox{ если } a \vee_{(o)} b = o\\
	~[ox], & \mbox{ если } a \vee_{(o)} b = x \ne o
	\end{array}
	\right.
\]
It follows from the definition of $\R$-tree, that in the first case
\[\rho(a,b) = \rho(a, o) + \rho(o, b),\]
and in the second
\[\rho(a,b) = \rho(a, x) + \rho(x, b).\]

\textit{Uniqueness}. Let $\preceq$ be upper rooted order on $X$ satisfying conditions of the theorem with root $o$. Assume that $x \preceq y$. Then it follows from (1) in definition \ref{dist} that $y \in [ox]$, and consequently $x \preceq_{(o)} y$. Conversely, let $x \preceq_{(o)} y$. Denote $x \vee y = w$ supremum of points $x$ and $y$ in the sense of the order $\preceq$. Then equalities hold:
\[\rho(x, w) + \rho(w,o) = \rho(x, o),\]
\[\rho(y,w) +\rho(w, o) = \rho(y,o)\]
and
\[\rho(x, w) + \rho(y, w) = \rho(x,y).\]
Since $X$ is $\R$-tree, we conclude that $w = y$ and $x \preceq y$.

\begin{remark}
In conditions of the theorem the semilattice $(X, \preceq_{(o)})$ can be lattice iff $X$ is one-point set, segment or semiinterval and $o$ is its end point. $X$ is linearly ordered in these situations.
\end{remark}

\textit{Sufficiency}. Let $X$ be a geodesic metric space endowed with upper-semilinear partial order $\preceq$ such that the triple $(X, \rho, \preceq)$ is a metric $\vee$-semilattice. Consider arbitrary points $x, y \in X$, for which $x \preceq y$. We have $y = x \vee y$. Let $z \in [xy]$ and assume that one of two relations $x \preceq z$ or $z \preceq y$ does not hold.

If the first relation does not valid and the second one does, then
\[x \preceq x \vee z \preceq y\]
and
\begin{align*}
&\rho(x,z) + \rho(z,y) = \rho(x, x \vee z + 2 \rho(z, x\vee z) + \rho(x \vee z, y) >\\ &> \rho(x, x \vee z) + \rho(x \vee z, y) \ge \rho(x,y).
\end{align*}

If the second relation does not hold, then
\[x \preceq y \preceq x \vee z\]
and
\[\rho(x, z) > \rho(x,y).\]

Inequalities in both cases contradict to equality
\[\rho(x,z) + \rho(z,y) = \rho(x,y),\]
necessary for inclusion $z \in [xy]$. Consequently for all points $z$ of segment $[xy]$ 
\[x \preceq z \preceq y\]
holds.

Let now points $x$ and $y$ be incomparable in relation $\preceq$. Then they are connected by segment consisting of pair of subsegments $[xw]$ and $[wy]$, where $w = x \vee y$. We show that the segment is unique connecting $x$ and $y$. Consider arbitrary point $z \in X$. It is true or false condition
\begin{equation}\label{alt}
z \preceq w.
\end{equation}
for it. If the condition \eqref{alt} is false, then
\[\rho(x,z) > \rho(x, w)\]
and
\[\rho(z,y) > \rho(w,y).\]
After addition we obtain
\[\rho(x,z) + \rho(z,y) > \rho(x,w) + \rho(w,y) = \rho(x,y),\]
that is the point $z$ belongs no segment ended in $x$ and $y$. Let the condition \eqref{alt} is true. Then
\begin{align*}
&\rho(x,z) + \rho(z,y) \ge \rho(x, x \vee z) + \rho(x \vee z, w) +\\ &+ \rho(w, z \vee y) + \rho(z \vee y, y) = \rho(x,y),
\end{align*}
and equality holds only in the case when the point $z$ belongs one of segments $[xw]$ or $[wy]$. The uniqueness of segment $[xy] = [xw] \cup [wy]$ follows from here.

We prove that for arbitrary different points $x,y,z \in X$ the segment $[xy]$ is contained in the union of segments $[xz] \cup [zy]$. Consider possible (up to notation) cases.

\textbf{1.} $x \prec z \prec y$. In this case $[xy] = [xz] \cup [zy]$.

\textbf{2.} Points $x$ and $z$ are incomparable and $x \prec w \prec y$, where $w =  x\vee z$. Then
\[ [xy] = [xw] \cup [wy] \subset [xz] \cup [zy],\]
\[ [xz] = [xw] \cup [wz] \subset [xy] \cup [zy]\]
and
\[ [zy] = [zw] \cup [wy] \subset [zx] \cup [zy].\]

\textbf{3.} Points $x$ and $z$ are incomparable and $x \prec y \preceq  x \vee z = y \vee z$. Then $[xz] = [xy] \cup [yz]$.

\textbf{4.} All points $x$, $y$ are $z$ mutually incomparable. We denote $v = x \vee y$ and $w = x \vee z$. Because of upper semilinearity of the order $\preceq$ points $v$ and $w$ are comparable. We assume that $v \preceq w$. Then, as it is easy to see, $w = y \vee z$. Consequently,
\[ [xy] = [xw] \cup [vw] \cup [wy] \subset [xz] \cup [zy],\]
\[ [xz] = [xw] \cup [wz] \subset [xy] \cup [yz]\]
and
\[ [zy] = [zw] \cup [vw] \cup [wy] \subset [xz] \cup [xy].\]
Cases 1--4 conclude all possible situations up to renotations.
\qed

Following example shows essentiality of condition of upper semilinearity in the theorem proved above.
\begin{example}
We define following partial order in the coordinate plane $A^2$ with coordinates $(x,y)$: we set $(x_1, y_1) \preceq (x_2, y_2)$ iff $x_1 \le x_2$ and $y_1 \le y_2$. It is easy to see that the pair $(A^2, \preceq)$ is $\vee$-semilattice (and even lattice). Also we consider the metric $\rho$ in $A^2$, generated by the norm $\|(x,y)\| = |x| + |y|$. 
The triple $(A^2, \rho, \preceq)$ is metric $\vee$-semilattice, the metric space $(A^2, \rho)$ is geodesic, but is not $\R$-tree. The order $\preceq$ is not semilinear.
\end{example}

\section{The topological space $\mathcal O_+(X)$}\label{spaceoforders}

Here we study the structure of the space $\mathcal O_+(X)$ of partial orders on complete locally compact $\R$-tree $X$ which define on $X$ upper semilinear metric $\vee$-semilattices. We denote $\mathcal O_+^r$ its subspace consisting of rooted orders. In fact, main results of this paragraph can be extended to general case of $\R$-tree, but one can not use the notion of metric compactification: it is well-defined only for locally compact complete metric spaces.

The partial order $\tau$ as binary relation is a closed subset of the metric square $X \times X$. The sum metric on $X \times X$ is of view $d_+$:
\[d_+((x_1, x_2), (y_1, y_2)) = \rho(x_1, y_1)+\rho^2(x_2,y_2),\]
and it generates Hausdorff distance on the family of closed subsets.
\begin{lemma}\label{metric}
For two rooted orders $\sigma, \tau \in \mathcal O_+^r$ the Hausdorff distance $\Hd(\sigma, \tau)$ is finite and equal to distance between their roots.
\end{lemma}
\proof If $x$ is the root of the order $\tau$ and $y$ --- the root of $\sigma$, then
\begin{equation}\label{hdrho}
\Hd(\tau, \sigma) \le \rho(x,y).
\end{equation}
Really, if  for points $s,t \in X$ we have $s \tau t$ but not $s \sigma t$, then the point $t$ lies in the segment $[xs]$ but not in $[ys]$. The distance from $t$ to the segment $[ys]$ is not greater then $\rho(t, s)$. Choose in $[ys]$ the point $w$ nearest to $t$. Since $s \sigma w$ and $d_+((s,t), (s,w)) = \rho(t,w) \le \rho(s,t)$, then for any $\varepsilon > 0$
\[(s,t) \in \mathcal N_{\rho(x,y)+\varepsilon}(\sigma).\]
Analogously, if $s \sigma t$ but not $s \tau t$, then for any $\varepsilon > 0$
\[(s,t) \in \mathcal N_{\rho(x,y)+\varepsilon}(\tau).\]

From the other hand roots satisfy relations $x \sigma y$ and $y \tau x$. If $\varepsilon < \rho(x,y)$, then the metric ball $B((y,x),\varepsilon)$ in sense of metric $d_+$ does not intersect $\sigma \subset X \times X$: for $(p,q) \in B((y,x), \varepsilon)$ the sum of distances $\rho(y,p) + \rho(q,x) \le \varepsilon$. Consequently
\begin{equation}\label{rhohd}
\rho(x, y) \le \Hd(\tau, \sigma).
\end{equation}
Inequalities \eqref{hdrho} and \eqref{rhohd} lead to expected equality. \qed

So, the subset $\mathcal O_+^r(X) \subset \mathcal O_+(X)$ with Hausdorff metric generated by the sum metric on $X \times X$ is isometric to $X$. Hence, $\mathcal O_+^r$ is also $\R$-tree. For arbitrary orders in $\mathcal O_+(X)$ the Hausdorff distance can be infinite, so Hausdorff metric does not generate uniquely defined topology on  $\mathcal O_+(X)$.
We will make $\mathcal O_+(X)$ the topological space, defining the base $\mathcal B$ of the topology as following. For a pair of different points $x, y \in X$ we set
 \[\mathcal U_{(x,y)} = \{\tau \in \mathcal O_+(X)\ |\ x \tau y\} \setminus \{\preceq_{(y)}\}.\]
The base $\mathcal B$ consists of all open sets of the metric space $\mathcal O_+^r$ and  various sets of type $\mathcal U_{(x,y)}$ for $x \ne y \in X$ (here we ignore the fact that coincidence of sets $\mathcal U_{(x,y)} = \mathcal U_{(x',y)}$ for $x \ne x'$ is possible). The fact that the family $\mathcal B$ generates a topology on $\mathcal O_+(X)$ is based on following lemma.
\begin{lemma}
For arbitrary sets $\mathcal U_1, \mathcal U_2 \in \mathcal B$ and order $\tau \in \mathcal U_1 \cap \mathcal U_2$ there exists such a set $\mathcal U \in \mathcal B$, that
\[\tau \in \mathcal U \subset \mathcal U_1 \cap \mathcal U_2.\]
\end{lemma}
\proof The statement of the lemma is evident if $\mathcal U_1$ and $\mathcal U_2$ both are open in $\mathcal O_+^r$. Since the intersection of set $\mathcal U_{(x,y)}$ with $\mathcal O_+^r$ is also open in $\mathcal O_+^r$, the statement is valid for any rooted order $\preceq_{(x)} \in \mathcal U_1 \cap \mathcal U_2 \cap \mathcal O_+^r$ as well. Let $\tau \in \mathcal U_1 \cap \mathcal U_2$ be non-rooted order. Assuming that $\mathcal U_1 = \mathcal U_{(x_1, y_1)}$ and $\mathcal U_2 = \mathcal U_{(x_2, y_2)}$, denote $y$ arbitrary point such that $y_1 \vee_\tau y_2 <_\tau y$ (such point exists, since $\tau$ is non-rooted order). Then
\[\tau \in \mathcal U_{(x_1, y)} = \mathcal U_{(x_2, y)} \subset \mathcal U_{(x_1, y_1)} \cap \mathcal U_{(x_2, y_2)},\]
and the conclusion of lemma follows. \qed

\begin{lemma}\label{Haus}
The topological space $\mathcal O_+(X)$ is Hausdorff.
\end{lemma}
\proof Subspace $\mathcal O_+^r(X)\subset \mathcal O_+(X)$ is Hausdorff as metric one. Let $\tau = \preceq_x \in \mathcal O_+^r(X)$ and $\sigma \in \mathcal O_+(X) \setminus \mathcal O_+^r(X)$. Choose a point $y \in X$ such that $x \sigma y$. Sets $\mathcal U_{(x,y)}$ and metric ball $B(\tau, \rho(x,y))$ in Hausdorff metric are non-intersecting neighbourhoods of orders $\tau$ and $\sigma$. Finally, let both orders $\tau$ and $\sigma$ are non-rooted. Choose a pair of points $x, y \in X$, for which $x \tau y$ but not $x \sigma y$. Denote $w = x \vee_\sigma y$. Sets $\mathcal U_1 = \mathcal U_{(w,y)}$ and $\mathcal U_2 = \mathcal U_{(y,w)}$ are non-intersecting neighbourhoods of orders $\tau$ and $\sigma$ correspondingly. \qed

Every $\R$-tree is $CAT(0)$-space, that is simply connected space, non-positively curved in the sense of A.D. Aleksandrov. The theory of Aleksandrov spaces is well-developed now (cf. \cite{BH} for example). For complete locally compact $CAT(0)$-spaces the metric compactification coincides with so-called \textit{geodesic compactification} with  \textit{cone topology}. Two rays $c, d : [0, +\infty) \to X$ are called \textit{asymptotic}, if the function $\rho(c(t), d(t))$ is bounded for $t \in [0, +\infty)$. This means that Hausdorff distance between rays $c$ and $d$ is finite:
\[\Hd(c,d) < +\infty.\]
The asymptoticity is equivalence on the set of rays in $X$. \textit{The geodesic boundary} $\partial_gX$ is defined as a set of equivalency classes in this relation. \textit{The cone topology} on $\overline{X}_g = X \cup \partial_gX$ is the topology of uniform convergence at bounded domains of natural parameterization of segments and rays. Precisely, the sequence of points $\{\xi_n\} \subset \overline{X}_g$ converges in cone topology to the point $\eta \in \overline{X}_g$ iff the sequence of natural parameterizations of segments (rays) $[o\xi_n]$ beginning in the marked point $o \in X$ converges to the natural parameterization of segment (ray) $[o\eta]$ uniformly at common bounded domains of parameters. 

The coincidence of compactifications means that the identity map $\Id: X \to X$ can be continued to homeomorphism $\overline{X}_m \to \overline{X}_g$. From now we will not tell the difference between metric and geodesic compactifications. We use notation $\overline{X}= \overline{X}_m = \overline{X}_g$. The boundary of the space $X$, i.e. its set of points at infinity $\overline{X}\setminus X$ is denoted as $\partial_\infty X$. In particular, every horofunction $\beta$ on $X$ is \textit{Busemann function}, defined from some ray $c: \R_+ \to X$ by equality
\[\beta(y) = \lim_{t \to +\infty}(\rho(y, c(t) - t).\]

Here we state homeomorphism $\Phi$ of topological space $\mathcal O_+(X)$ to compactification $\overline{X}$.

Let a partial order $\tau$ on $\R$-tree $X$ defines upper semilinear metric $\vee$-semilattice. If $\tau$ is upper rooted order with root $x_\tau$, then we set $\Phi(\tau) = x_\tau$. If $\tau$ is not rooted, then its upper semilinearity and metric property imply, that for arbitrary point $x \in X$ the set $U_x$ is a union of mutually included into each other segments of type $[xy]$. Since $X$ is assumed to be complete and the order $\tau$ is non-rooted, the lengths of segments $[xy]$ increases to infinity. Such a union is a ray in $X$ beginning in $x$. For any two points $x, y \in X$ the supremum $w = x \vee y$ is defined and rays $U_x$ and $U_y$ intersect by the ray $U_w$, i.e. are asymptotic. In this case we set $\Phi(\tau) = [U_x]\in \partial_\infty X$, i.e. equivalence class of rays in $X$, which are asymptotic to $U_x$.

\begin{theorem}
The map $\Phi:\mathcal O_+ \to \overline{X}$ is a homeomorphism.
\end{theorem}

The map $\Phi$ is injective. Really, if $\tau_1,\tau_2 \in \mathcal O_+(X)$ are two different rooted orders on $X$, then from the uniqueness statement in theorem \ref{crit} their roots are different. Let orders $\tau_1, \tau_2$ be non-rooted. Consider a pair $x,y \in X$, such that $ x \tau_1 y$ but not $x \tau_2 y$. The set $\overline{X} \setminus \{y\}$ is subdivided to linear connectedness components, so that the point $x$ and the set $\{z \ |\ y \tau_1 z\}  \setminus \{y\}$ belong to different components. The image $\Phi(\tau_1)$ belongs to the closure of the set $\{z \ |\ y \tau_1 z\}$ in $\overline{X}$, while  $\Phi(\tau_2)$ to the closure of the component containing $x$.

The map $\Phi$ is surjective. Really, for arbitrary point $x \in X$, the theorem \ref{crit} defines the order $\tau \in \mathcal O_+$, such that $\Phi(\tau) = x$. For arbitrary point $\xi \in \partial X$ the order $\tau$ is defined within corresponding Busemann function $\beta_\xi$. For a point $y \in X$, \textit{the horoball}, that is the sublevel set
\[\hb(\xi, y) = \{z \in X\ |\ \beta_\xi(z) \le \beta_\xi(y)\}\]
corresponding to it, is a convex set in $X$, and for $x \in X$ its projection  $\pi_{\xi, y}(x)$ to $\hb(\xi, y)$, that is the nearest to $x$ point of $\hb(\xi, y)$ is uniquely defined. We set $x \tau_\xi y$ iff $\pi_{\xi, y}(x) = y$. The relation $\tau_\xi$ is an order in $\mathcal O_+(X)$ and $\Phi(\tau_\xi) = \xi$.

So the map $\Phi$ is bijective. According to the lemma \ref{Haus} and compactness of $\overline{X}$ it is sufficient to prove that  $\Phi$ is open map for proving the theorem. In view of lemma \ref{metric} and definition of the topology on $\mathcal O_+(X)$, it is sufficient to verify that for any pair $x, y \in X$ the image of the set $\mathcal U_{(x,y)}$ is open in $\overline{X}$. We prove that the complement $\overline{X} \setminus \Phi(\mathcal U_{(x,y)}$ is closed.

For $\tau \in \mathcal U_{(x,y)}\cap \mathcal O_+^r$ denote $z$ the root of $\tau$ and $r = \rho(y,z)$. If for a point $t \in X$ the inequality $\rho(z, t) < r$ holds, then $y \in [xt]$ and $x \preceq_t y$. Consequently, $t \in \mathcal U_{(x,y)}$, and $\tau$ is interior point of the set $\mathcal U_{(x,y)}$. So, $X \setminus \mathcal U_{(x,y)}$ is closed in $X$ and contains all limit points, belonging to $X$. 

Let the sequence $\{\xi_n\} \subset \overline{X} \setminus \mathcal U_{(x,y)}$ converges to the ideal point $\xi \in \partial X$. We may assume that the marked point $o \notin \Phi(\mathcal U_{(x,y)})$. In such assumption for all $n$ the equality $[o\xi_n] \cap \Phi(\mathcal U_{(x,y)}) = \varnothing$ holds. Let $R > \rho(o,y)$ and $c_n: [0, R] \to X$ be natural parameterizations of segments or rays $[o\xi_n]$. Then points $c_n(R)$ converge to the point $c(R)$, where $c: \R_+ \to X$ is natural parameterization of the ray $[o\xi]$. This remains true if we choose arbitrary number $R' > R$. But since $X$ is $\R$-tree, then it is possible only in the case when $c_n(R) = c(R)$ for all $n$, greater then some sufficiently large number $N \in \N$. Consequently, $c(R) \notin \Phi(\mathcal U_{(x,y)})$. Since, in particular, $c(R) \in [y \xi]$, the ray $[y \xi]$ is asymptotic to no ray with beginning segment $[xy]$. We conclude that $\xi \notin \Phi(\mathcal U_{(x,y)})$ and the set $\overline{X} \setminus \Phi(\mathcal U_{(x,y)})$ is closed.

\section{Similarity homogeneous non-homogeneous $\R$-tree}\label{example}

In this paragraph we construct the series of the metric spaces $X$, which are similarity homogeneous non-homogeneous $\R$-trees. This proves essentiality of the condition of locally compactness in Berestovski\v\i's conjecture. Our examples has closed relation with universal $\R$-trees $A_\mu$ (here $\mu$ is a cardinal), constructed by A.Dyubina and I.Polterovich in \cite{DP}. Our space is conformally equivalent to the universal $\R$-trees $A_\mu$ in the sense of definition of conform equivalence from \cite{Be}.

The function $y = f(x)$, defined in interval $(\alpha, \beta)$ is called \textit{piecewise constant from the left}, if for every point $x\in (\alpha, \beta)$ there exists $\varepsilon > 0$, such that $f|_{[x-\varepsilon, x]} = \const$.
The set $X$ is a set of pairs $(f, a_f)$, where $a_f > 0$ is real number and $f: (a_f, +\infty) \to G$ is piecewise constant from the left function with target an additive group $G$ such that $f|_{(b_f, +\infty)} \equiv 0$ for some $b_f \ge a_f$.

We introduce the partial order on $X$ putting $(f, a_f) \preceq (g, a_g)$ iff $a_f \le a_g$ and $f|_{(a_g, +\infty)} = g$.

\begin{lemma}
The pair $(X, \preceq)$ is upper semilinear $\vee$-semilattice.
\end{lemma}

\proof For pairs $(f, a_f),(g, a_g) \in X$ their supremum is a pair $(h, a_h)$, where
\[a_h = \inf\{x \in \R\ |\ f|_{(x, +\infty)} = g|_{(x, +\infty)}\}\]
and $h = f|_{(a_h, +\infty)} = g|_{(a_h, +\infty)}$.
For pair $(f, a_f)$ linearity of the order of upper cone $U_{(f, a_f)}$ follows from the definition of  $\preceq$. \qed

Next we define a metric $\rho$ on the set $X$, putting
\[\rho((f, a_f), (g, a_g)) = |a_f - a_g|,\]
if $(f, a_f) \preceq (g, a_g)$ or $(g, a_g) \preceq (f, a_f)$, and
\[\rho(p,q) = \rho(p, p \vee q) + \rho(p \vee q, q),\]
if elements $p,q \in X$ incomparable.

\begin{lemma}
The metric space $(X, \rho)$ is geodesic. A triple $(X, \rho, \preceq)$ is upper metric semilattice.
\end{lemma}

\proof It is sufficient to show that if $(f, a_f) \preceq (g, a_g)$, then they can be connected by a segment. Such segment is represented by the parameterization $\gamma: [a_f, a_g] \to X$ given by formula
\[\gamma(t) = (f|_{(t, +\infty)}, t).\]
Conditions (1) and (2) in the definition \ref{dist} follows from representation of the metric $\rho$.
\qed

\begin{corollary}
The space $X$ is an $\R$-tree.
\end{corollary}

The point $x$ in $\R$-tree $X$ is called \textit{point of valency $\mu$}, where $\mu$ is a cardinal, if the cardinal number of connected components of $X \setminus \{x\}$ is $\mu$.

\begin{lemma}
Every point of $\R$-tree $X$ is a point of valency $|G|+1$.
\end{lemma}

\proof For the point $x = (f, a_f)$ connected components of the set $X \setminus \{x\}$ are: 
\[A = \{(g, a_g)\ |\ \mbox{ не выполняется } (g, a_g) \prec (f, a_f)\}\]
and for every element $\alpha \in G$
\[B_\alpha = \{(g, a_g)\ |\ (g, a_g)\prec (f, a_f) \mbox{ и } g(a_f) = \alpha\}.\]
Totally $|G|+1$ components. \qed

\begin{lemma}\label{submetr}
The map $\varphi: X \to \R$ given by formula $\varphi((f, a_f) = a_f$ is a submetry.
\end{lemma}

\proof It follows from the definition of the metric $\rho$. \qed

\begin{lemma}
$\R$-tree $X$ is locally complete. The completeness radius in the point $(f, a_f)\in X$ is  $r_c(f, a_f) = a_f$.
\end{lemma}

\proof Let $\{(f_n, a_{f_n})\}$ be a fundamental sequence of points of the ball $B((f, a_f), r)$ where $r < a_f$. This means that for any $\varepsilon > 0$ there exists such a number $N \in \N_1$ that for all $m,n > N_1$ functions $f_m$ and $f_n$ coincide in the ray $[a_{f_m}+\frac13 \varepsilon, +\infty)$. From lemma \ref{submetr} the sequence $\{a_{f_n}\}$ is also fundamental. Let $a_g = \lim\limits_{n \to \infty} a_n$. Now define the function $g: (a_g, +\infty) \to \{0,1\}$. For arbitrary $t > a_g$ take $\varepsilon = t - a_g$ and choose such $N_2 \in \N$ that $a_{f_n} - a_g < \frac13 \varepsilon$ for all $n > N_2$. Then for $n > N = \max\{N_1, N_2\}$ functions $f_n(t)$ stabilize: $f_n(t) = f_{N+1}(t)$, at that there exists $\varepsilon_1 < \frac13\varepsilon$, such that $f_n|_{[t-\varepsilon_1]} = f_n(t)$. Putting $g(t) = \lim\limits_{n \to \infty} f_n(t)$, we obtain well-defined piecewise constant from the left function $g$ on the ray $(a_g, + \infty)$ for which $\rho((g, a_g), (f_n, a_{f_n}))$ converges to zero when $n \to \infty$. It follows that the pair $(g,a_g)$ is an element of the space $X$ and the limit of the sequence $\{(f_n, a_{f_n})\}$. Moreover, from above it immediately follows the estimation on the completeness radius $r_c(f, a_f) \ge a_f$. The inverse inequality is obvious. \qed

Next we consider following transformations of $\R$-tree $X$.

For a number $\lambda > 0$ put $H_\lambda(f, a_f) = (f_\lambda, \lambda a_f)$, where $f_\lambda(x) = f(\frac1\lambda x)$. The transformation $H_\lambda$ is a similarity with  coefficient $\lambda$, that is, equality $\rho(H_\lambda(x), H_\lambda(y)) = \lambda \rho(x,y)$ holds for all elements $x,y \in X$.

Let us extend an arbitrary function $f: (a_f, +\infty) \to G$ to the function $\bar f: \R_+ \to G$, putting $\bar f(t) = 0$ for $t \le a_f$. The pair $(f, a_f)$ generates the transformation $R_f$ of the space $X$ by formula
\[R_f((g, a_g)) = (g + \bar f|_{(a_g, +\infty)}, a_g).\]
Automatically $R_f$ is an isometry and for the function $\mathbf{0}$, identically equal to zero one have $R_f((\mathbf{0}, a_f)) = (f, a_f)$. Transformations $R_{f_1}$ and $R_{f_2}$ coincide when and only when $\bar f_1 = \bar f_2$.
We note that $R_{f_2}\circ R_{f_1} = R_{f_1}\circ R_{f_2} = R_{f_1+ f_2}$ и $(R_f)^{-1} = R_{-f}$, hence the set $\mathcal R$ of isometries $R_f$ gives commutative subgroup of $\Isom(X)$.

\begin{theorem}
$X$ is similarity homogeneous non-homogeneous $\R$-tree. It is a metric fibration with fibers
\[F_{a} = \{(f, a)\}.\]
\end{theorem}

\proof For points $(f, a_f)$ and $(g, a_g)$ the superposition $R_g \circ H_\lambda \circ (R_f)^{-1}$ where $\lambda = a_g/a_f$, is a similarity moving $(f, a_f)$ в $(g, a_g)$.

Sets $F_a$ are closed because their complements are open. The map $T_{(a, b)} : F_{a} \to F_{b}$ defined by
\[T_{(a,b)}(f,a) = (T_{a-b}f, b),\]
where the function $T_{a-b}f: (b, +\infty) \to G$ acts by the formula $T_{a-b}f(x) = f(x+a-b)$, is an isometry of the fiber $F_a$ onto the fiber $F_b$. We assume that $a < b$. For the point $(f, a) \in F_a$ its nearest point in the fiber $F_b$ is the point $(f|_{(b, +\infty)}, b)$. Such nearest point is unique. For the point $(f, b)\in F_b$ nearest points of the fiber $F_a$ are points of type $(g, a)$ under the condition $g|_{(b, +\infty)} = f$. So all conditions of the definition of metric fibration are fulfilled. \qed

The theorem 2.1 of the paper \cite{Be} states that locally complete similarity homogeneous metric space is homogeneous if and only if it is complete. The completeness radius in the point $x = (f, a_f)\in X$ is equal to $c(x)=a_f$ so $X$ is not homogeneous. The space  $X$ is not homeomorphic to product $F_a \times \R_+$: it is linearly complete, while the factor $F_a$ in the product $F_a \times \R_+$, is not connected because it is ultrametric in the metric induced from  $X$.

Splitting \eqref{split} of the group $\Sim(X)$ is generated by inclusion homomorphism  $\Isom(X) \to \Sim(X)$ and projection homomorphism $\Sim(X) \to \R_+$ comparing to every similarity its coefficient. Such splitting does not generate the structure of topological product $\Sim(X) = \Isom(X)\times \R_+$: subgroup in $\Sim(X)$ generated all transformations of type $H_\lambda$ and $R_f$ is linearly connected, while $\Isom(X)$ acts on ultrametric spaces $F_a$ as totally non-connected group. Every connected component in the product $\Isom(X) \times \R_+$ is isomorphic to $\R_+$.

\end{document}